\documentclass[12pt,oneside,reqno]{amsart}
\usepackage{mathrsfs}
\usepackage{amssymb}
\newcommand{\dif}{\mathrm{d}}

\newcommand{\be}{\begin{eqnarray}}
\newcommand{\ee}{\end{eqnarray}}
\newcommand{\ce}{\begin{eqnarray*}}
\newcommand{\de}{\end{eqnarray*}}
\newtheorem{theorem}{Theorem}[section]
\newtheorem{lemma}[theorem]{Lemma}
\newtheorem{remark}[theorem]{Remark}
\newtheorem{definition}[theorem]{Definition}
\newtheorem{proposition}[theorem]{Proposition}
\newtheorem{Examples}[theorem]{Examples}
\newtheorem{corollary}[theorem]{Corollary}

\def\[{{\Big[}}
\def\]{{\Big]}}
\def\<{{\langle}}
\def\>{{\rangle}}
\def\({{\Big(}}
\def\){{\Big)}}

\def\tr{{\rm tr}}

\def\no{\nonumber}
\def\bt{\begin{theorem}}
\def\et{\end{theorem}}
\def\bl{\begin{lemma}}
\def\el{\end{lemma}}
\def\br{\begin{remark}}
\def\er{\end{remark}}
\def\bx{\begin{Examples}}
\def\ex{\end{Examples}}
\def\bd{\begin{definition}}
\def\ed{\end{definition}}
\def\bp{\begin{proposition}}
\def\ep{\end{proposition}}
\def\bc{\begin{corollary}}
\def\ec{\end{corollary}}

\def\cF{{\mathcal F}}

\def\cU{{\mathcal U}}

\def\mE{{\mathbb E}}

\def\mN{{\mathbb N}}

\def\mR{{\mathbb R}}

\def\mU{{\mathbb U}}

\def\geq{\geqslant}
\def\leq{\leqslant}

\begin{document}

\allowdisplaybreaks

\title{Exponential ergodicity for SDEs with jumps and non-Lipschitz coefficients*}

\author{Huijie Qiao}



\thanks{*This work is supported by NSF(No. 11001051) of China.}


\date{}

\dedicatory{Department of Mathematics,
Southeast University,\\
Nanjing, Jiangsu 211189, P.R.China\\
Email: hjqiaogean@yahoo.com.cn}

\begin{abstract}
In this paper we show irreducibility and the strong Feller property
for transition probabilities of stochastic differential equations
with jumps and monotone coefficients. Thus, exponential ergodicity
and the spectral gap for the corresponding transition semigroups are
obtained.
\end{abstract}

\maketitle \rm

\section{Introduction}
Let $(\Omega,\cF,P;(\cF_t)_{t\geq 0})$ be a complete filtered
probability space, and $(\mU,\cU,\nu)$ a $\sigma$-finite measurable
space. Let $\{W(t)\}_{t\geq 0}$ be a $d$-dimensional standard
$\cF_t$-adapted Brownian motion, and $\{k_t,t\geq 0\}$ a stationary
$\cF_t$-adapted Poisson point process with values in $\mU$ and with
characteristic measure $\nu$(cf. \cite{iw}). Let $N_{k}((0,t],\dif
u)$ be the counting measure of $k_{t}$, i.e., for $A\in\cU$
$$
N_k((0,t],A):=\#\{0<s\leq t: k_s\in A\},
$$
where $\#$ denotes the cardinality of a set. The compensator measure
of $N_k$ is given by
$$
\tilde{N}_{k}((0,t],\dif u):=N_{k}((0,t],\dif u)-t\nu(\dif u).
$$

In the following, we fix a $\mU_{0}\in\cU$ such that
$\nu(\mU-\mU_{0})<\infty$, and consider the following stochastic
differential equation (SDE) with jumps in $\mR^{d}$: \be
X_t=x_0+\int^t_0b(X_s)\dif s+\int^t_0\sigma(X_s)\dif W_s
+\int^{t+}_0\int_{\mU_{0}}f(X_{s-},u)\tilde{N}_{k}(\dif s\dif u),
\label{Eqj1} \ee where $b:\mR^d\to\mR^d$,
$\sigma:\mR^d\to\mR^{d\times d}$ and $f:\mR^d\times\mU\to\mR^d$ are
measurable functions. Here, the second integral of the right side in
Eq.(\ref{Eqj1}) is taken in the It\^o's sense, and the definition of
the third integral is referred to \cite{iw}. When $d=1$, Youngmff
Kwon and Chaniio I. Ff in \cite{kw} showed that the transition
semigroup of Eq.(\ref{Eqj1}) is strong Feller and irreducible under
some smoothness and growth conditions on $b, \sigma, f$ with
nondegenerate diffusion term. If $b, \sigma, f$ are Lipschitz
continuous, Masuda in \cite{ma1} provided sets of conditions under
which the transition semigroup of Eq.(\ref{Eqj1}) fulfils the
ergodic theorem for any initial distribution.

In this paper, we study the ergodicity of Eq.(\ref{Eqj1}) under some
non-Lipschitz conditions. First of all, recall some notions about
the ergodicity. Let $\{X_t(x),t\geq0,x\in\mE\}$ be a family of
Markov processes with state space $\mE$ being a Hausdorff topology
space, and
transition probability $p_t(x,E)$. Then\\
(i) $p_t$ is called irreducible if for each $t>0$ and $x\in\mE$ \ce
p_t(x,E)>0~\mbox{ for any non-empty open set $E\subset\mE$ }; \de
(ii) $p_t$ is called strong Feller if for each $t>0$ and
$E\in\mathscr{B}(\mE)$ \ce \mE\ni x\mapsto p_t(x,E)\in[0,1]~\mbox{
is continuous }; \de (iii) A measure $\mu$ on
$(\mE,\mathscr{B}(\mE))$ is an invariant measure for $p_t$ if \ce
\int_{\mE}p_t(x,E)\mu(\dif x)=\mu(E), \quad\forall t>0,
E\in\mathscr{B}(\mE). \de

The transition probability $p_t(x,\cdot)$ determines a Markov
semigroup $(p_t)_{t\geq 0}$. The theorem below is a classical result
combining the above concepts.(cf. \cite{dz})

\bt\label{fii} Assume a Makov semigroup $(p_t)_{t\geq 0}$ is
irreducible and strong Feller. Then there exists at most one
invariant measure for it. Moreover, if $\mu$ is the invariant
measure, then $\mu$ is ergodic and equivalent to each $p_t(x,\cdot)$
and as $t\rightarrow\infty$, $p_t(x,E)\rightarrow\mu(E)$ for any
Borel set $E$. \et
Next introduce our non-Lipschitz conditions.\\
{\bf Hypotheses:}

({\bf H}$_1$) There exists $\lambda_0\in\mR$ such that for all $x,y\in\mR^d$
\ce
2\<x-y, b(x)-b(y)\>+\|\sigma(x)-\sigma(y)\|^2\leq
\lambda_0|x-y|^2\kappa(|x-y|),
\de
where $\kappa$ is a positive continuous function, bounded on $[1,\infty)$ and satisfying
\ce
\lim\limits_{x\downarrow0}\frac{\kappa(x)}{\log
x^{-1}}=\delta<\infty.
\de

Here the function $\kappa$ controls the continuity modulus of $b(x)$
and $\sigma(x)$ such that the modulus is non-Lipschitz, for example,
$\kappa(x)= C_1\cdot(\log(1/x)\vee K)^{1/\beta_1}$ for some
$\beta_1>1$ and $C_1,K>0$.

({\bf H}$_2$) There exists $\lambda_1>0$
such that for all $x\in\mR^d$
$$
|b(x)|^2+\|\sigma(x)\|^2\leq \lambda_1(1+|x|)^2.
$$

({\bf H}$_3$) $b$ is continuous and there exists $\lambda_2>0$ such that
\be
\<\sigma(x)h,h\>\geq\sqrt{\lambda_2}|h|^2, \qquad
x, h\in\mR^d. \label{sig}
\ee

Here $\<\cdot,\cdot\>$ denotes the inner product in $\mR^d$, $|\cdot|$
the length of a vector in $\mR^d$ and $\|\cdot\|$ the Hilbert-Schmit
norm from $\mR^d$ to $\mR^d$.

When $\nu(\mU_{0})=0$ and $b, \sigma$ satisfy the above assumptions
and another assumption, Zhang \cite{z1} proved that the transition
semigroup of Eq.(\ref{Eqj1}) has the exponential ergodicity in the sense
that there exists a constant $\beta_2>0$ such that for $t>0$,
\ce
\|p_t(x_0,\cdot)-\mu\|_{Var}\leq C_2\cdot e^{-\beta_2t},
\de
where $\|\cdot\|_{Var}$ denotes the total variation of a signed measure
and $C_2>0$ is a constant. Here we
require $\nu(\mU_{0})\neq0$ and follow the same lines as done in
\cite{z1}. Thus how to treat the term with jumps is our key.

Firstly, we make the following assumption on $f$:

({\bf H$_f$}) \ce \int_{\mU_0}\big|f(x,u)-f(y,u)\big|^2\nu(\dif
u)\leq2|\lambda_0||x-y|^2\kappa(|x-y|) \de and for $q=2$ and $4$
$$
\int_{\mU_{0}}|f(x,u)|^{q}\,\nu(\dif u)\leq \lambda_1(1+|x|)^{q}.
$$
Under ({\bf H}$_1$), ({\bf H}$_2$) and {\bf (H$_f$)}, it is well
known that there exists a unique strong solution to Eq.(\ref{Eqj1})
(cf. \cite[Theorem 170, p.140]{si}). This solution will be denoted
by $X_t(x_0)$. The transition semigroup associated with $X_t(x_0)$
is defined by \ce p_t\varphi(x_0):={\bf E}\varphi(X_t(x_0)), \qquad
t>0, \quad\varphi\in B_b(\mR^d), \de where $B_b(\mR^d)$ stands for
the Banach space of all bounded measurable functions on $\mR^d$. The
transition probability is given by \ce
p_t(x_0,E):=(p_t1_E)(x_0)=P(X_t(x_0)\in E), \quad
E\in\mathscr{B}(\mR^d). \de By Girsanov's theorem on processes with
jumps we get the irreducibility under ({\bf H}$_1$)-({\bf H}$_3$)
and {\bf (H$_f$)} (cf. Proposition \ref{ir}).

In order to get strong Feller property, we need the following
stronger assumptions on $b, \sigma, f$:

({\bf H}$_1^\prime$) There exists $\lambda_0\in\mR$ such that for
all $x,y\in\mR^d$ \ce 2\<x-y,
b(x)-b(y)\>+\|\sigma_{\lambda_2}(x)-\sigma_{\lambda_2}(y)\|^2\leq
\lambda_0|x-y|^2\kappa(|x-y|), \de where $\sigma_{\lambda_2}(x)$ is
the unique symmetric nonnegative definite matrix-valued function
such that $\sigma_{\lambda_2}(x)\sigma_{\lambda_2}(x)
=\sigma(x)\sigma(x)-\lambda_2\emph{I}$ for the unit matrix $\emph{I}$.

({\bf H$_f^\prime$}) There exists a positive function $L(u)$
satisfying \ce \sup_{u\in\mU_0}L(u)\leq\gamma<1~\mbox{ and }
\int_{\mU_0}L(u)^2\,\nu(\dif u)<+\infty, \de such that for any
$x,y\in\mR^{d}$ and $u\in\mU_0$ \ce |f(x,u)-f(y,u)|\leq L(u)|x-y|,
\de and \ce |f(0,u)|\leq L(u). \de

\br\label{sigma} To explain that ({\bf H}$_1^\prime$) is stronger
than ({\bf H}$_1$), a matrix result is needed. And we will give a
general result in Section 2. \er

When $b,\sigma$ and $f$ satisfy ({\bf H}$_1^\prime$), ({\bf H}$_2$),
({\bf H}$_3$) and {\bf (H$_f^\prime$)}, we obtain strong Feller
property by the coupling method (cf. Proposition \ref{sf}).

Finally, to show the ergodicity,
the following assumption is needed:\\
({\bf H$_{b,\sigma,f}$}) There exist a $r\geq 2$ and two constants
$\lambda_3>0$, $\lambda_4\geq0$ such that for all $x\in\mR^d$ \ce
2\<x,b(x)\>+\|\sigma(x)\|^2+\int_{\mU_0}\big|f(x,u)\big|^2\nu(\dif
u)\leq-\lambda_3|x|^r+\lambda_4. \de We are now in a position to state our main
result in the present paper. \bt\label{er} Assume ({\bf
H}$_1^\prime$), ({\bf H}$_2$), ({\bf H}$_3$) and {\bf
(H$_f^\prime$)}. Then the semigroup $p_t$ is irreducible and strong
Feller. If in addition, ({\bf H}$_{b,\sigma,f}$)  hold, then there
exists a unique invariant probability measure $\mu$ of $p_t$
having full support in $\mR^d$ such that

(i) if $r\geq2$ in ({\bf H}$_{b,\sigma,f}$), then for all $t>0$ and
$x_0\in\mR^d$, $\mu$ is equivalent to $p_t(x_0,\cdot)$ and
\ce
\lim\limits_{t\rightarrow\infty}\|p_t(x_0,\cdot)-\mu\|_{Var}=0.
\de

(ii) if $r>2$ in ({\bf H}$_{b,\sigma,f}$), then for some $\alpha,
C_3>0$ independent of $x_0$ and $t$,
\ce
\|p_t(x_0,\cdot)-\mu\|_{Var}\leq C_3\cdot e^{-\alpha t}.
\de
Moreover,
for any $\gamma>1$ and each $\varphi\in L^\gamma(\mR^d,\mu)$
\ce
\|p_t\varphi-\mu(\varphi)\|_{L^\gamma(\mR^d,\mu)}\leq C_4\cdot
e^{-\alpha t/\gamma}\|\varphi\|_{L^{\gamma}(\mR^d,\mu)},
\quad\forall t>0,
\de
where $\alpha$ is the same as above, $\mu(\varphi):=\int_{\mR^d}\varphi(x)\mu(\dif x)$
and $C_4>0$ is a constant depending on $\gamma$. In particular,
let $L_{\gamma}$ be the generator of $p_t$ in
$L^{\gamma}(\mR^d,\mu)$, then $L_{\gamma}$ has a spectral gap
(greater than $\alpha/\gamma$) in $L^{\gamma}(\mR^d,\mu)$. \et

The following convention will be used throughout the paper: $C$ with
or without indices will denote different positive constants
(depending on the indices) whose values may change from one place to
another one.

\section{Proof of Theorem \ref{er}}

\bl\label{mr} Suppose $A,B$ are two symmetric positive definite
matrices satisfying that there exists a $\lambda>0$ such that
$\<Ah,h\>\geq\sqrt{\lambda}|h|^2, \<Bh,h\>\geq\sqrt{\lambda}|h|^2$
for $h\in\mR^d$, and $AB=BA$. Then
\be
\|A-B\|\leq\|\sqrt{A^2-\lambda\emph{I}}-\sqrt{B^2-\lambda\emph{I}}\|.
\label{norine}
\ee
\el
\begin{proof}
Set
\ce
A_\lambda:=\sqrt{A^2-\lambda\emph{I}}, \qquad B_\lambda:=\sqrt{B^2-\lambda\emph{I}}.
\de
To show (\ref{norine}), we consider the difference of $\|A-B\|^2$ and
$\|A_\lambda-B_\lambda\|^2$, i.e.
\be
\|A-B\|^2-\|A_\lambda-B_\lambda\|^2&=&\tr((A-B)^2)-\tr((A_\lambda-B_\lambda)^2)\no\\
&=&\tr(A^2-2AB+B^2)-\tr(A_\lambda^2-2A_\lambda B_\lambda+B_\lambda^2)\no\\
&=&\tr(A^2-2AB+B^2)-\tr(A^2-\lambda\emph{I}-2A_\lambda B_\lambda+B^2-\lambda\emph{I})\no\\
&=&\tr(-2AB+2\lambda\emph{I}+2A_\lambda B_\lambda)\no\\
&=&2\left[\tr(A_\lambda B_\lambda)-\tr(AB)+\lambda d\right],
\label{difest}
\ee
where $\tr(\cdot)$ stands for the trace of a matrix.

By the proof of \cite[Theorem 7.4.10, p.433]{hj}, one can obtain
$A=UMU^*$ and $B=UNU^*$, where $U\in\mR^{d\times d}$ is real orthogonal,
$M=diag(\eta_1,\dots,\eta_d)$, $N=diag(\mu_1,\dots,\mu_d)$, and $\eta_i, \mu_i$
are eigenvalues of $A$ and $B$, respectively, and larger than $\sqrt{\lambda}$.
Moreover, $A_\lambda=UM_\lambda U^*$
and $B_\lambda=UN_\lambda U^*$, where $M_\lambda=diag(\sqrt{\eta_1^2-\lambda},\dots,
\sqrt{\eta_d^2-\lambda})$, $N_\lambda=diag(\sqrt{\mu_1^2-\lambda},\dots,\sqrt{\mu_d^2-\lambda})$.
Thus,
\ce
\tr(A_\lambda B_\lambda)=\tr(UM_\lambda U^*UN_\lambda U^*)=\tr(UM_\lambda N_\lambda U^*)
=\sum\limits_{i=1}^d\sqrt{\eta_i^2-\lambda}\sqrt{\mu_i^2-\lambda}.
\de
By the same deduction as the above one, we have $\tr(AB)=\sum\limits_{i=1}^d\eta_i\mu_i$.
So, the right hand side of (\ref{difest}) can be written as $2\sum\limits_{i=1}^d
\big(\sqrt{\eta_i^2-\lambda}\sqrt{\mu_i^2-\lambda}-(\eta_i\mu_i-\lambda)\big)$.

Noting that
\ce
\left(\sqrt{\eta_i^2-\lambda}\sqrt{\mu_i^2-\lambda}\right)^2-(\eta_i\mu_i-\lambda)^2
&=&(\eta_i^2-\lambda)(\mu_i^2-\lambda)-(\eta_i\mu_i-\lambda)^2\\
&=&\eta_i^2\mu_i^2-\eta_i^2\lambda-\lambda\mu_i^2+\lambda^2-\eta_i^2\mu_i^2
+2\lambda\eta_i\mu_i-\lambda^2\\
&=&-\lambda(\eta_i-\mu_i)^2\leq0,
\de
we get $\sqrt{\eta_i^2-\lambda}\sqrt{\mu_i^2-\lambda}\leq\eta_i\mu_i-\lambda$. Thus,
it holds that $\|A-B\|^2-\|A_\lambda-B_\lambda\|^2\leq0$ and (\ref{norine}).
\end{proof}

\vspace{3mm}

To show the irreducibility, we firstly estimate $X_t(x_0)$.

\bl\label{cb} Suppose that $b,\sigma$ and $f$ satisfy ({\bf H}$_1$)-({\bf H}$_2$) and
{\bf (H$_f$)}. Then for any $T>0$
\be
{\bf E}\left[\sup\limits_{t\in[0,T]}|X_t|^2\right]\leq C,
\label{l2}
\ee
where $C$ depends on $x_0,\lambda_1$ and $T$.
\el
\begin{proof}
Applying It\^o's formula to Eq.(\ref{Eqj1}), we have
\be
|X_t|^2&=&|x_0|^2+\int_0^t2\<X_s,b(X_s)\>\dif s+\int_0^t2\<X_s,\sigma(X_s)\dif W_s\>\no\\
&&+\int_0^t\int_{\mU_0}\left[|X_{s-}+f(X_{s-},u)|^2-|X_{s-}|^2\right]\tilde{N}_k(\dif s,\dif u)+\int_0^t\|\sigma(X_s)\|^2\dif s\no\\
&&+\int_0^t\int_{\mU_0}\left[|X_s+f(X_{s},u)|^2-|X_s|^2-2\<X_s,f(X_{s},u)\>\right]\nu(\dif
u)\dif s.
\label{ito1}
\ee
By BDG inequality, Young's inequality and mean theorem, one can obtain that
\be
{\bf E}\left[\sup\limits_{t\in[0,T]}|X_t|^2\right]
&\leq&|x_0|^2+{\bf E}\int_0^T2|X_s||b(X_s)|\dif s
+C{\bf E}\left(\int_0^T|X_s|^2\|\sigma(X_s)\|^2\dif s\right)^{\frac{1}{2}}\no\\
&&+C{\bf E}\left[\int_0^T\int_{\mU_0}\left[|f(X_{s-},u)|^2+|f(X_{s-},u)||X_{s-}|\right]^2N_k(\dif s,\dif u)\right]^{\frac{1}{2}}\no\\
&&+{\bf E}\int_0^T\|\sigma(X_s)\|^2\dif s+{\bf E}\int_0^T\int_{\mU_0}|f(X_{s},u)|^2\nu(\dif u)\dif s\no\\
&\leq&|x_0|^2+{\bf E}\left(2\sup\limits_{t\in[0,T]}|X_t|\int_0^T|b(X_s)|\dif s\right)
+{\bf E}\int_0^T\|\sigma(X_s)\|^2\dif s\no\\
&&+C{\bf E}\left(\sup\limits_{t\in[0,T]}|X_t|^2\int_0^T\|\sigma(X_s)\|^2\dif s\right)^{\frac{1}{2}}
+{\bf E}\int_0^T\int_{\mU_0}|f(X_{s},u)|^2\nu(\dif u)\dif s\no\\
&&+C{\bf E}\left[\int_0^T\int_{\mU_0}2\left[|f(X_{s-},u)|^4+|f(X_{s-},u)|^2|X_{s-}|^2\right]N_k(\dif s,\dif u)\right]^{\frac{1}{2}}\no\\
&\leq&|x_0|^2+\frac{1}{4}{\bf E}\left[\sup\limits_{t\in[0,T]}|X_t|^2\right]
+C{\bf E}\int_0^T|b(X_s)|^2\dif s+C{\bf E}\int_0^T\|\sigma(X_s)\|^2\dif s\no\\
&&+C{\bf E}\left[\int_0^T\int_{\mU_0}|f(X_{s-},u)|^4N_k(\dif s,\dif u)\right]^{\frac{1}{2}}
+{\bf E}\int_0^T\int_{\mU_0}|f(X_{s},u)|^2\nu(\dif u)\dif s\no\\
&&+C{\bf E}\left[\int_0^T\int_{\mU_0}|f(X_{s-},u)|^2|X_{s-}|^2N_k(\dif s,\dif u)\right]^{\frac{1}{2}}.
\label{e11}
\ee

For the fifth term in the right hand side, we use
Young's inequality and {\bf (H$_f$)} to get
\be
&&C{\bf E}\left[\int_0^T\int_{\mU_0}|f(X_{s-},u)|^4N_k(\dif s,\dif u)\right]^{\frac{1}{2}}\no\\
&\leq&C{\bf E}\left[\sup\limits_{t\in[0,T]}(1+|X_{t-}|)^2\int_0^T\int_{\mU_0}\frac{|f(X_{s-},u)|^4}{(1+|X_{s-}|)^2}N_k(\dif s,\dif u)\right]^{\frac{1}{2}}\no\\
&\leq&\frac{1}{8}{\bf E}\left[\sup\limits_{t\in[0,T]}(1+|X_{t-}|)^2\right]+C{\bf E}\int_0^T\int_{\mU_0}\frac{|f(X_{s-},u)|^4}{(1+|X_{s-}|)^2}\nu(\dif u)\dif s\no\\
&\leq&\frac{1}{4}{\bf E}\left[\sup\limits_{t\in[0,T]}|X_{s-}|^2\right]+\frac{1}{4}+C{\bf
E}\int_0^T(1+|X_s|)^2\dif s.
\label{e12}
\ee

To the seventh term in
the right hand side, the similar method yields that \be
&&C{\bf E}\left[\int_0^T\int_{\mU_0}|f(X_{s-},u)|^2|X_{s-}|^2N_k(\dif s,\dif u)\right]^{\frac{1}{2}}\no\\
&\leq&\frac{1}{4}{\bf E}\left[\sup\limits_{t\in[0,T]}|X_{t-}|^2\right]+C{\bf
E}\int_0^T\int_{\mU_0}|f(X_s,u)|^2\nu(\dif u)\dif s.
\label{e3}
\ee

Combining (\ref{e11}), (\ref{e12}), (\ref{e3}), ({\bf H}$_2$) and
{\bf (H$_f$)}, we obtain \ce
{\bf E}\left[\sup\limits_{t\in[0,T]}|X_t|^2\right]&\leq&4|x_0|^2+1+C{\bf E}\int_0^T(1+|X_s|)^2\dif s\\
&\leq&4|x_0|^2+1+C{\bf E}\int_0^T2(1+|X_s|^2)\dif s\\
&\leq&4|x_0|^2+(1+2CT)+2C\int_0^T{\bf E}\left[\sup\limits_{s\in[0,t]}|X_s|^2\right]\dif t.
\de
Gronwall's inequality yields
\ce
{\bf E}\left[\sup\limits_{t\in[0,T]}|X_t|^2\right]\leq C,
\de
where $C$ depends on $x_0, \lambda_1$ and $T$.
\end{proof}

Secondly, construct some auxiliary processes. For any $T>0$, let $t_0\in(0,T)$, whose value will
be determined below. Set for any $n\in\mN$
$$
X_{t_0}^n=X_{t_0}I_{\{|X_{t_0}|\leq n\}}.
$$
And then by Lemma \ref{cb}
$$
\lim\limits_{n\rightarrow\infty}{\bf E}|X_{t_0}^n-X_{t_0}|^2=0.
$$
For $t\in[t_0,T]$ and $y\in\mR^d$, define
\ce
&&J_t^n=\frac{T-t}{T-t_0}X_{t_0}^n+\frac{t-t_0}{T-t_0}y,\\
&&h_t^n=\frac{y-X_{t_0}^n}{T-t_0}-b(J_t^n).
\de
Thus,
\ce
J_{t_0}^n=X_{t_0}^n, \quad J_T^n=y,
\de
and $J_t^n$ satisfies the following equation:
\ce
J_t^n=X_{t_0}^n+\int_{t_0}^tb(J_s^n)\dif s+\int_{t_0}^th_s^n\dif s, \quad t\in[t_0,T].
\de

Next, we introduce the following equation:
\ce
Y_t&=&X_{t_0}+\int_{t_0}^tb(Y_s)\dif s+\int_{t_0}^th_s^n\dif s+\int_{t_0}^t\sigma(Y_s)\dif W_s\\
&&+\int^{t+}_{t_0}\int_{\mU_{0}}f(Y_{s-},u)\,\tilde{N}_{k}(\dif s\dif u), \qquad t\in[t_0,T].
\de

\bl\label{con} Suppose that $b,\sigma$ and $f$ satisfy ({\bf H}$_1$)-({\bf H}$_2$) and
{\bf (H$_f$)}. Then
\be
{\bf E}|Y_T-y|^2\leq\left[{\bf E}|X_{t_0}-X_{t_0}^n|^2+C(T-t_0)\right]^{e^{-|\lambda_0|(T-t_0)}}.
\label{l1}
\ee
\el
\begin{proof}
Set $Z_t:=Y_t-J_t^n$, and then $Z_t$ satisfies the following equation
\ce
Z_t&=&X_{t_0}-X_{t_0}^n+\int_{t_0}^t\left(b(Y_s)-b(J_s^n)\right)\dif s+\int_{t_0}^t\sigma(Y_s)\dif W_s\\
&&+\int^{t+}_{t_0}\int_{\mU_{0}}f(Y_{s-},u)\,\tilde{N}_{k}(\dif s\dif u), \qquad t\in[t_0,T].
\de
By It\^o's formula we obtain
\ce
|Z_t|^2&=&|X_{t_0}-X_{t_0}^n|^2+\int_{t_0}^t2\<Z_s,b(Y_s)-b(J_s^n)\>\dif s+\int_{t_0}^t2\<Z_s,\sigma(Y_s)\dif W_s\>\no\\
&&+\int_{t_0}^t\int_{\mU_0}\left[|Z_{s-}+f(Y_{s-},u)|^2-|Z_{s-}|^2\right]\tilde{N}_k(\dif s,\dif u)+\int_{t_0}^t\|\sigma(Y_s)\|^2\dif s\no\\
&&+\int_{t_0}^t\int_{\mU_0}\left[|Z_s+f(Y_{s},u)|^2-|Z_s|^2-2\<Z_s,f(Y_{s},u)\>\right]\nu(\dif
u)\dif s.
\de
It follows from ({\bf H}$_1$)-({\bf H}$_2$) and {\bf (H$_f$)} that
\be
{\bf E}|Z_t|^2&=&{\bf E}|X_{t_0}-X_{t_0}^n|^2+{\bf E}\int_{t_0}^t2\<Z_s,b(Y_s)-b(J_s^n)\>\dif s+{\bf E}\int_{t_0}^t\|\sigma(Y_s)\|^2\dif s\no\\
&&+{\bf E}\int_{t_0}^t\int_{\mU_0}\left[|Z_s+f(Y_{s},u)|^2-|Z_s|^2-2\<Z_s,f(Y_{s},u)\>\right]\nu(\dif u)\dif s\no\\
&\leq&{\bf E}|X_{t_0}-X_{t_0}^n|^2+{\bf E}\int_{t_0}^t|\lambda_0||Z_s|^2\kappa(|Z_s|)\dif s+2\lambda_1{\bf E}\int_{t_0}^t(1+|Y_s|)^2\dif s.
\label{zong}
\ee

There exists a $\delta>0$ such that
\be
x^2\kappa(x)\leq\rho_\delta(x^2), \qquad x>0,
\label{non}
\ee
where $\rho_\delta: \mR_+\mapsto\mR_+$ is a concave function given by
\ce \rho_\delta(x):=\left\{
\begin{array}{lcl}
x\log x^{-1},&& x\leq\delta,\\
(\log\eta^{-1}-1)x+\eta, &&x>\delta.
\end{array}
\right.
\de

Next, estimate ${\bf E}|Y_t|^2$. For $t\in[t_0,T]$, by H\"older's inequality, one can get
\ce
|Y_t|^2&\leq&5|X_{t_0}|^2+5T^{\frac{1}{2}}\int_{t_0}^t|b(Y_s)|^2\dif s+5T^{\frac{1}{2}}\int_{t_0}^t|h_s^n|^2\dif s
+5\left|\int_{t_0}^t\sigma(Y_s)\dif W_s\right|^2\\
&&+5\left|\int^{t+}_{t_0}\int_{\mU_{0}}f(Y_{s-},u)\,\tilde{N}_{k}(\dif s\dif u)\right|^2.
\de
Moreover, by Burkholder's inequality, ({\bf H}$_2$) and {\bf (H$_f$)}, it holds that
\ce
{\bf E}|Y_t|^2&\leq&5{\bf E}|X_{t_0}|^2+5T^{\frac{1}{2}}{\bf E}\int_{t_0}^t|b(Y_s)|^2\dif s+5T^{\frac{1}{2}}{\bf E}\int_{t_0}^t|h_s^n|^2\dif s\\
&&+5\int_{t_0}^t\|\sigma(Y_s)\|^2\dif s+5\int^{t+}_{t_0}\int_{\mU_{0}}|f(Y_{s-},u)|^2\nu(\dif u)\dif s\\
&\leq&5{\bf E}|X_{t_0}|^2+C+C\int_{t_0}^t{\bf E}|Y_s|^2\dif s.
\de
Gronwall's inequality admits us to obtain that
\be
\sup\limits_{t\in[t_0,T]}{\bf E}|Y_t|^2\leq C,
\label{ye}
\ee
where $C$ is independent of $t_0$.

Combining (\ref{zong})-(\ref{ye}), by Jensen's inequality and the Bihari inequality
(cf. \cite[Lemma 2.1]{z1}), we have
\ce
{\bf E}|Y_T-y|^2\leq\left[{\bf E}|X_{t_0}-X_{t_0}^n|^2+C(T-t_0)\right]^{e^{-|\lambda_0|(T-t_0)}}.
\de
The proof is completed.
\end{proof}

\bp\label{ir} Suppose that $b,\sigma$ and $f$ satisfy ({\bf
H}$_1$)-({\bf H}$_3$) and {\bf (H$_f$)}. Then the transition
probability $p_t$ is irreducible. \ep
\begin{proof}
To prove the irreducibility, it suffices to prove that for each
$T>0$ and $x_0\in\mR^d$,
\ce
p_T\big(x_0,B(y,a)\big)=P\big(X_T(x_0)\in B(y,a)\big)=P\big(|X_T(x_0)-y|<a\big)>0,
\de
or equivalently
\be
P\big(|X_T(x_0)-y|\geq a\big)<1,
\label{irre}
\ee
for any $y\in\mR^d$ and $a>0$.

First of all, study the process $Y_t$. Define
\ce
Y_t:=X_t, \qquad t\in[0,t_0]
\de
and then
\ce
Y_t&=&x_0+\int_0^tb(Y_s)\dif s+\int_0^tI_{\{s>t_0\}}h_s^n\dif s+\int_0^t\sigma(Y_s)\dif W_s\\
&&+\int^{t+}_{t_0}\int_{\mU_{0}}f(Y_{s-},u)\,\tilde{N}_{k}(\dif s\dif u), \qquad t\in[0,T].
\de

Set \ce
&&H_t:=I_{\{t>t_0\}}\sigma(Y_t)^{-1}h_t^n,\\
&&\xi_t:=\exp\left\{\int_0^{t}\<\dif
W_s,H_s\>-\frac{1}{2}\int_0^{t} |H_s|^2\dif s\right\}.
\de By ({\bf H}$_3$), we obtain that $|H_{t}|^2$ is
bounded, which yields that ${\bf E}\xi_T=1$ by Novikov's criteria.
And then define \ce
&&\bar{W}_t:=W_t+\int_0^{t}H_s\dif s,\\
&&Q:=\xi_TP. \de Thus by \cite[Theorem 132]{si} we know that $Q$ is
a probability measure, $\bar{W}_t$ is a $Q$ -Brownian motion and
$\tilde{N}_k\big((0,t],\dif u\big)$ is a Poisson martingale measure
under $Q$ with the same compensator $\nu(\dif u)t$. Moreover, $Y_t$
is the solution of the following equation
\ce
Y_{t}=x_0+\int^{t}_{0}b(Y_{s})\,\dif s+\int^{t}_{0}\sigma(Y_{s})\,\dif \bar{W}_{s}
+\int^{t+}_{0}\int_{\mU_{0}}f(Y_{s-},u)\,\tilde{N}_{k}(\dif s\dif u).
\de

By the uniqueness in law of Eq.(\ref{Eqj1}) we attain that the law of $\{X_t, t\in[0,T]\}$
under $P$ is the same to that of $\{Y_t, t\in[0,T]\}$ under $Q$. Therefore, to obtain
(\ref{irre}), it is sufficient to prove $Q\big(|Y_T-y|\geq a\big)<1$, and furthermore, $P\big(|Y_T-y|\geq a\big)<1$
by equivalency of $Q, P$.

It holds by Chebyshev's inequality and Lemma \ref{con} that
\ce
P\big(|Y_T-y|\geq a\big)\leq\left[{\bf E}|X_{t_0}-X_{t_0}^n|^2+C(T-t_0)\right]^{e^{-|\lambda_0|(T-t_0)}}/a^2.
\de
Choosing $n$ large enough and $t_0$ close enough to $T$, we have
$$
P\big(|Y_T-y|\geq a\big)<1.
$$
The proof is completed.
\end{proof}

\vspace{3mm}

Next we use the coupling method to prove strong Feller property. Set
$a(x):=\sigma(x)\sigma(x)$. And then the infinitesimal generator of
Eq.(\ref{Eqj1}) is given by \ce
L\psi(x)&=&b^i(x)\partial_i\psi(x)+\frac{1}{2}a^{ij}(x)\partial_{ij}\psi(x)\\
&&+\int_{\mU_0}\left(\psi\left(x+f(x,u)\right)-\psi(x)-f^i(x,u)\partial_i\psi(x)\right)\nu(\dif
u), \de for $\psi\in C_b^2(\mR^d)$ and $x\in\mR^d$. Recall that an
operator $\tilde{L}$ on $\mR^{2d}$ is called a coupling operator of
$L$ if $\tilde{L}$ satisfies the marginal condition: \be
(\tilde{L}\psi)(x,y)=L\psi(x), \quad \psi\in C_b^2(\mR^d), \quad
x,y\in\mR^d, \label{cp} \ee where $\psi$ is regarded as a function in $C_b^2(\mR^{2d})$. For
any $\delta\in(0,1)$ and $|x_0-y_0|<\delta$, we define \ce
u(x,y):=\frac{x-y}{|x-y|}, \quad
u_\delta(x,y)=\frac{|x_0-y_0|^{\frac{\alpha}{2}}}
{\delta^{\frac{\alpha}{2}}}u(x,y),\\
c(x,y):=\lambda_2\big(I-2u_\delta(x,y)u_\delta(x,y)^*\big)
+\sigma_{\lambda_2}(x)\sigma_{\lambda_2}(y)^*, \de for $x,y\in\mR^d,
x\neq y$, where $\alpha\in(0,1)$ (its value will be determined
below). Thus the operator given by \ce
\tilde{L}\psi(x,y)&=&b^i(x)\partial_{x_i}\psi(x,y)+b^i(x)\partial_{y_i}\psi(x,y)\\
&&+\frac{1}{2}\[a^{ij}(x)\partial_{x_ix_j}\psi(x,y)
+a^{ij}(x)\partial_{y_iy_j}\psi(x,y)\\
&&+c^{ij}(x,y)\partial_{x_iy_j}\psi(x,y)+(c(x,y)^*)^{ij}\partial_{y_ix_j}\psi(x,y)\]\\
&&+\int_{\mU_0}\(\psi\big(x+f(x,u),y+f(y,u)\big)-\psi(x,y)-f^i(x,u)\partial_{x_i}\psi(x,y)\\
&&\qquad -f^i(y,u)\partial_{y_i}\psi(x,y)\)\nu(\dif u), \quad
\psi\in C_b^2(\mR^{2d}) \de is a coupling operator of $L$.

By the analysis similar to \cite[Section 3.1]{wa}, there exist a
stochastic basis
$(\tilde{\Omega},\tilde{\cF},\tilde{P};(\tilde{\cF}_t)_{t\geq 0})$
and a $\hat{\mR}^{2d}$-valued process $\tilde{Z}^\prime_t$
($\hat{\mR}^{2d}$ is one point compactification of $\mR^{2d}$) such
that $\tilde{Z}^\prime_0=(x_0,y_0)$ and for each $\psi\in
C_b^2(\mR^{2d})$ \be \left\{\psi(\tilde{Z}^\prime_{t\wedge
e})-\psi(\tilde{Z}^\prime_0)-\int_0^{t\wedge e}\tilde{L}
\psi(\tilde{Z}^\prime_{s})\dif s, \quad t\geq 0\right\} \label{lm}
\ee is an $(\tilde{\cF}_t)$-local martingale, where $e$ is the
explosion time of the process
$\tilde{Z}^\prime_t=(\tilde{X}^\prime_t,\tilde{Y}_t)$, i.e. \ce
e=\lim\limits_{n\rightarrow\infty}e_n, \quad e_n=\inf\left\{t>0:
\left(|\tilde{X}^\prime_t|+|\tilde{Y}_t|\right)>n\right\}. \de Since
$\tilde{X}^\prime_t$ and $\tilde{Y}_t$ are associated with the
operator $L$ starting from $x_0$ and $y_0$, respectively, we have
$e=\infty$ by ({\bf H}$_2$) and ({\bf H$_f^\prime$}). Consider the
coupling time $\tau$ of $\tilde{Z}^\prime_t$, i.e.
$$
\tau:=\inf\{t>0:|\tilde{X}^\prime_t-\tilde{Y}_t|=0\},
$$
and define \ce
\tilde{X}_t=\left\{\begin{array}{l}\tilde{X}^\prime_t, \quad t<\tau,\\
\tilde{Y}_t, ~\quad t\geq\tau.
\end{array}
\right. \de And then for any $\psi\in C^\infty_c(\mR^d)$, \ce
\psi(\tilde{X}_t)-\psi(x_0)-\int_0^tL\psi(\tilde{X}_s)\dif s
&=&\psi(\tilde{X}_{t\wedge\tau})-\psi(x_0)-\int_0^{t\wedge\tau}L\psi(\tilde{X}_s)\dif s\\
&&+\psi(\tilde{X}_t)-\psi(\tilde{X}_{t\wedge\tau})-\int_{t\wedge\tau}^tL\psi(\tilde{X}_s)\dif s\\
&=&\psi(\tilde{X}^\prime_{t\wedge\tau})-\psi(x_0)-\int_0^{t\wedge\tau}L\psi(\tilde{X}^\prime_s)\dif s\\
&&+\psi(\tilde{Y}_t)-\psi(\tilde{Y}_{t\wedge\tau})-\int_{t\wedge\tau}^tL\psi(\tilde{Y}_s)\dif
s, \de together with Doob stopping theorem yields that the process
$\tilde{X}_t$ is also associated with $L$. Thus the generator of
$\tilde{Z}_t:=(\tilde{X}_t,\tilde{Y}_t)$ before $\tau$ is just
$\tilde{L}$.

\bp\label{sf} Under ({\bf H}$_1^\prime$), ({\bf H}$_2$), ({\bf
H}$_3$) and {\bf (H$_f^\prime$)}, the semigroup $p_t$ of $X_t$ is
strong Feller. \ep
\begin{proof}
For any $\varphi\in B_b(\mR^d)$, by the definition of $\tilde{X}_t$
\ce
|p_t\varphi(x_0)-p_t\varphi(y_0)|&=&|\tilde{{\bf E}}(\varphi(\tilde{X}_t(x_0)))-\tilde{{\bf E}}(\varphi(\tilde{Y}_t(y_0)))|\\
&\leq&|\tilde{{\bf E}}[\varphi(\tilde{X}_t(x_0))1_{t<\tau}]-\tilde{{\bf E}}[\varphi(\tilde{Y}_t(y_0))1_{t<\tau}]|\\
&&+|\tilde{{\bf E}}[\varphi(\tilde{X}_t(x_0))1_{t\geq\tau}]-\tilde{{\bf E}}[\varphi(\tilde{Y}_t(y_0))1_{t\geq\tau}]|\\
&\leq&2\|\varphi\|_0\tilde{P}(t<\tau), \de where $\tilde{{\bf E}}$
is the expectation with respect to $\tilde{P}$.

Next we estimate $\tilde{P}(t<\tau)$.

Define \ce S_\delta:=\inf\{t\geq0:
|\tilde{X}_t-\tilde{Y}_t|>\delta\}. \de Setting
$g(r):=\frac{r}{1+r}, r\geq0$ and $\psi(x,y):=g(|x-y|)$, one can
obtain by (\ref{lm}) \be
&&\tilde{{\bf E}}g(|\tilde{X}_{t\wedge\tau\wedge e_n\wedge S_\delta}-\tilde{Y}_{t\wedge\tau\wedge e_n\wedge S_\delta}|)\no\\
&=&g(|x_0-y_0|)+\tilde{{\bf E}}\int_0^{t\wedge\tau\wedge e_n\wedge
S_\delta}\frac{\bar{G}(\tilde{X}_s,\tilde{Y}_s)}{2}
g^{\prime\prime}(|\tilde{X}_s-\tilde{Y}_s|)\dif s\no\\
&&+\tilde{{\bf E}}\int_0^{t\wedge\tau\wedge e_n\wedge
S_\delta}\frac{\tr\left(G(\tilde{X}_s,\tilde{Y}_s)\right)
-\bar{G}(\tilde{X}_s,\tilde{Y}_s)+2F(\tilde{X}_s,\tilde{Y}_s)}{2|\tilde{X}_s-\tilde{Y}_s|}g^{\prime}(|\tilde{X}_s-\tilde{Y}_s|)\dif s\no\\
&&+\tilde{{\bf E}}\int_0^{t\wedge\tau\wedge e_n\wedge
S_\delta}\int_{\mU_0}\[g\left(|\tilde{X}_s+f(\tilde{X}_s,u)-\tilde{Y}_s-f(\tilde{Y}_s,u)|\right)
-g\left(|\tilde{X}_s-\tilde{Y}_s|\right)\no\\
&&\qquad-g^{\prime}\left(|\tilde{X}_s-\tilde{Y}_s|\right)\left<u(\tilde{X}_s,\tilde{Y}_s),
f(\tilde{X}_s,u)-f(\tilde{Y}_s,u)\right>\]\nu(\dif u)\dif s\no\\
&=:&g(|x_0-y_0|)+I_1+I_2+I_3, \label{ito} \ee where \ce
&&G(x,y)=a(x)+a(y)-c(x,y)-c(x,y)^*,\\
&&\bar{G}(x,y)=\<u(x,y),G(x,y)u(x,y)\>,\\
&&F(x,y)=\<x-y,b(x)-b(y)\>, \de and \ce
g^{\prime}(r)=\frac{1}{(1+r)^2}, \quad
g^{\prime\prime}(r)=-\frac{2}{(1+r)^3}. \de

Noting that by ({\bf H}$_2$) \ce
\bar{G}(x,y)&=&\big<u(x,y),\big(a(x)+a(y)-2\lambda_2
-\sigma_{\lambda_2}(x)\sigma_{\lambda_2}(y)-\sigma_{\lambda_2}(y)\sigma_{\lambda_2}(x)\big)u(x,y)\big>\\
&&+\<u(x,y),4\lambda_2\frac{|x_0-y_0|^\alpha}{\delta^\alpha}u(x,y)u(x,y)^*u(x,y)\>\\
&=&\left<u(x,y),\big(\sigma_{\lambda_2}(x)-\sigma_{\lambda_2}(y)\big)^2u(x,y)\right>
+4\lambda_2\frac{|x_0-y_0|^\alpha}{\delta^\alpha}\\
&\geq&4\lambda_2\frac{|x_0-y_0|^\alpha}{\delta^\alpha}, \de and \ce
\tr\big(G(x,y)\big)-\bar{G}(x,y)&=&\tr\big(a(x)+a(y)-2\lambda_2
-\sigma_{\lambda_2}(x)\sigma_{\lambda_2}(y)-\sigma_{\lambda_2}(y)\sigma_{\lambda_2}(x)\big)\\
&&+\tr\big(4\lambda_2\frac{|x_0-y_0|^\alpha}{\delta^\alpha}u(x,y)u(x,y)^*\big)
-\bar{G}(x,y)\\
&\leq&\|\sigma_{\lambda_2}(x)-\sigma_{\lambda_2}(y)\|^2. \de Thus
\be
I_1\leq-\frac{4\lambda_2}{\delta^\alpha(1+\delta)^3}\cdot|x_0-y_0|^\alpha\cdot
\tilde{{\bf E}}(t\wedge\tau\wedge e_n\wedge S_\delta). \label{i1}
\ee By ({\bf H}$_1^\prime$) we have \be I_2\leq\tilde{{\bf
E}}\int_0^{t\wedge\tau\wedge e_n\wedge
S_\delta}\frac{|\lambda_0|}{2}
|\tilde{X}_s-\tilde{Y}_s|\kappa(|\tilde{X}_s-\tilde{Y}_s|)\dif s.
\label{i2} \ee Next, deal with $I_3$. Mean theorem and ({\bf
H$_f^\prime$}) admit us to get \ce
&&g(|x-y+f(x,u)-f(y,u)|)-g(|x-y|)-g^\prime(|x-y|)\<u(x,y),f(x,u)-f(y,u)\>\\
&=&g^{\prime\prime}(|x-y+\theta(f(x,u)-f(y,u))|)|x-y+\theta(f(x,u)-f(y,u))|^{-2}\\
&&\quad \cdot\left<x-y+\theta(f(x,u)-f(y,u)),f(x,u)-f(y,u)\right>^2\\
&&+g^{\prime}(|x-y+\theta(f(x,u)-f(y,u))|)\[|x-y+\theta(f(x,u)-f(y,u))|^{-1}\\
&&\quad \cdot|f(x,u)-f(y,u)|^2-\<x-y+\theta(f(x,u)-f(y,u)),f(x,u)-f(y,u)\>^2\\
&&\quad \cdot|x-y+\theta(f(x,u)-f(y,u))|^{-3}\]\\
&\leq&\frac{L(u)^2}{1-L(u)}|x-y|\\
&\leq&C L(u)^2|x-y|, \de where $0<\theta<1$. So \be I_3\leq
\tilde{{\bf E}}\int_0^{t\wedge\tau\wedge e_n\wedge
S_\delta}C|\tilde{X}_s -\tilde{Y}_s|\dif s. \label{i3} \ee

Combining (\ref{i1}), (\ref{i2}) and (\ref{i3}), one can obtain that
\be &&\frac{1}{1+\delta}\tilde{{\bf E}}|\tilde{X}_{t\wedge\tau\wedge
e_n\wedge S_\delta}
-\tilde{Y}_{t\wedge\tau\wedge e_n\wedge S_\delta}|\no\\
&\leq&|x_0-y_0|-\frac{4\lambda_2}{\delta^\alpha(1+\delta)^3}\cdot|x_0-y_0|^\alpha\cdot
\tilde{{\bf E}}(t\wedge\tau\wedge e_n\wedge S_\delta)\no\\
&&+\frac{|\lambda_0|}{2}\int_0^{t\wedge\tau\wedge e_n\wedge
S_\delta}\tilde{{\bf E}} \rho_\delta(|\tilde{X}_s-\tilde{Y}_s|)\dif
s+C\int_0^{t\wedge\tau\wedge e_n\wedge S_\delta}
\tilde{{\bf E}}|\tilde{X}_s-\tilde{Y}_s|\dif s\no\\
&\leq&|x_0-y_0|+(\frac{|\lambda_0|}{2}+C)\int_0^{t\wedge\tau\wedge
e_n\wedge S_\delta}
\rho_\delta(\tilde{{\bf E}}|\tilde{X}_s-\tilde{Y}_s|)\dif s\no\\
&&-\frac{4\lambda_2}{\delta^\alpha(1+\delta)^3}\cdot|x_0-y_0|^\alpha\cdot
\tilde{{\bf E}}(t\wedge\tau\wedge e_n\wedge S_\delta),
\label{itoend} \ee
where the second step bases on Jensen's inequality. By the
Bihari inequality (cf. \cite[Lemma 2.1]{z1}) , we get that for any
$t>0$ and $|x_0-y_0|\leq\delta$, \be \tilde{{\bf
E}}|\tilde{X}_{t\wedge\tau\wedge e_n\wedge S_\delta}
-\tilde{Y}_{t\wedge\tau\wedge e_n\wedge S_\delta}|\leq
(1+\delta)|x_0-y_0|^{\exp\{-(1+\delta)(\frac{|\lambda_0|}{2}+C)t\}}.
\label{me} \ee

On one hand, substituting (\ref{me}) into (\ref{itoend}) yields \ce
\tilde{{\bf E}}(t\wedge\tau\wedge e_n\wedge S_\delta)
&\leq&\frac{\delta^\alpha(1+\delta)^3}{4\lambda_2}\[|x_0-y_0|^{1-\alpha}+(\frac{|\lambda_0|}{2}+C)t\\
&&\quad \cdot\rho_\delta
\left((1+\delta)|x_0-y_0|^{\exp\{-(1+\delta)(\frac{|\lambda_0|}{2}+C)t\}}\right)\cdot|x_0-y_0|^{-\alpha}\].
\de Letting $n\rightarrow\infty$, one can obtain by Levy's theorem
\ce \tilde{{\bf E}}(t\wedge\tau\wedge S_\delta)
&\leq&\frac{\delta^\alpha(1+\delta)^3}{4\lambda_2}\[|x_0-y_0|^{1-\alpha}+(\frac{|\lambda_0|}{2}+C)t\\
&&\quad \cdot\rho_\delta
\left((1+\delta)|x_0-y_0|^{\exp\{-(1+\delta)(\frac{|\lambda_0|}{2}+C)t\}}\right)\cdot|x_0-y_0|^{-\alpha}\].
\de Take $\alpha=\exp\{-(1+\delta)(|\lambda_0|+2C)t\}/3$. Thus,
there exists a $0<\delta'<\delta$ such that for any
$|x_0-y_0|\leq\delta'$ \be \tilde{{\bf E}}\big((2t)\wedge\tau\wedge
S_\delta\big)\leq
C_{t,\lambda_0,\delta}\cdot|x_0-y_0|^{\exp\{-(1+\delta)(|\lambda_0|+2C)t\}/6}.
\label{e1} \ee

On the other hand, it follows from (\ref{me}) \ce
\delta\tilde{P}\big(t\wedge\tau\wedge e_n>S_\delta\big)
&\leq&\tilde{{\bf E}}\[|\tilde{X}_{t\wedge\tau\wedge e_n\wedge
S_\delta} -\tilde{Y}_{t\wedge\tau\wedge e_n\wedge S_\delta}|
\cdot I_{t\wedge\tau\wedge e_n>S_\delta}\]\\
&\leq&\tilde{{\bf E}}|\tilde{X}_{t\wedge\tau\wedge e_n\wedge
S_\delta}
-\tilde{Y}_{t\wedge\tau\wedge e_n\wedge S_\delta}|\\
&\leq&(1+\delta)|x_0-y_0|^{\exp\{-(1+\delta)(\frac{|\lambda_0|}{2}+C)t\}}.
\de Letting $n\rightarrow\infty$, we have \be
\tilde{P}\big((2t)\wedge\tau>S_\delta\big)
\leq\frac{1+\delta}{\delta}|x_0-y_0|^{\exp\{-(1+\delta)(|\lambda_0|+2C)t\}}.
\label{e2} \ee

Finally, by (\ref{e1}) and (\ref{e2}) it holds that \ce
\tilde{P}(t<\tau)&=&\tilde{P}(t<\tau, S_\delta>t)+\tilde{P}(t<\tau, S_\delta\leq t)\\
&\leq&\tilde{P}\big((2t)\wedge\tau\wedge S_\delta>t\big)+\tilde{P}\big((2t)\wedge\tau>S_\delta\big)\\
&\leq&\frac{1}{t}\tilde{{\bf E}}\big((2t)\wedge\tau\wedge
S_\delta\big)
+\frac{1+\delta}{\delta}|x_0-y_0|^{\exp\{-(1+\delta)(|\lambda_0|+2C)t\}}\\
&\leq&C_{t,\lambda_0,\delta}\cdot|x_0-y_0|^{\exp\{-(1+\delta)(|\lambda_0|+2C)t\}/6}.
\de Thus, the proof is completed.
\end{proof}

\vspace{3mm}

We now give

\vspace{3mm}

\textbf{Proof of Theorem \ref{er}}

\vspace{3mm}

By It\^o's formula and mean theorem we get \be
{\bf E}|X_t|^2&=&|x_0|^2+\int_0^t{\bf E}2\<X_s,b(X_s)\>\dif s+\int_0^t{\bf E}\|\sigma(X_s)\|^2\dif s\no\\
&&+\int_0^t\int_{\mU_0}{\bf E}\left[|X_s+f(X_{s},u)|^2-|X_s|^2-2\<X_s,f(X_{s},u)\>\right]\nu(\dif u)\dif s\no\\
&=&|x_0|^2+\int_0^t{\bf E}2\<X_s,b(X_s)\>\dif s+\int_0^t{\bf E}\|\sigma(X_s)\|^2\dif s\no\\
&&+\int_0^t\int_{\mU_0}{\bf E}|f(X_{s},u)|^2\nu(\dif u)\dif s.
\label{ito3} \ee By ({\bf H}$_{b,\sigma,f}$), one obtains \ce
\frac{1}{t}\int_0^t{\bf E}|X_s|^r\dif
s\leq\frac{x_0}{\lambda_3t}+\frac{\lambda_4}{\lambda_3}. \de

Set
\ce
\mu_T(A):=\frac{1}{T}\int_0^Tp_t(x_0,A)\dif t,
\de
for any $T>0$ and $A\in\mathscr{B}(\mR^d)$. And we have by Chebyshev's inequality
\ce
\mu_T(B^c(0,R))&=&\frac{1}{T}\int_0^Tp_t(x_0,B^c(0,R))\dif t\\
&\leq&\frac{1}{TR^r}\int_0^T{\bf E}|X_t|^r\dif t\\
&\leq&\frac{1}{R^r}\left(\frac{x_0}{\lambda_3T}+\frac{\lambda_4}{\lambda_3}\right).
\de
Thus, for any $\varepsilon>0$, $\mu_T(B(0,R))>1-\varepsilon$ for
$R$ being large enough. Hence, $\{\mu_T,T>0\}$ is tight and its
limit is an invariant probability measure $\mu$. As we have just
proved, $p_t$ is irreducible and strong Feller, then by Theorem
\ref{fii}, $\mu$ is equivalent to each $p_t(x,\cdot)$ and as
$t\rightarrow\infty$, $p_t(x,E)\rightarrow\mu(E)$ for any Borel set
$E$.

If $r>2$, to (\ref{ito3}), taking derivatives
with respect to $t$ and using ({\bf H}$_{b,\sigma,f}$) and
H\"older's inequality give \ce
\frac{\dif{\bf E}|X_t|^2}{\dif t}&\leq&-\lambda_3{\bf E}|X_t|^r+\lambda_4\\
&\leq&-\lambda_3({\bf E}|X_t|^2)^{\frac{r}{2}}+\lambda_4. \de Let $f(t)$ solve
the following ODE: \ce \left\{
\begin{array}{lcl}
f'(t)=-\lambda_3f(t)^{\frac{r}{2}}+\lambda_4,\\
f(0)=|x_0|^2.
\end{array}
\right. \de By the comparison theorem on ODE and \cite[Lemma
1.2.6, p.32]{c}, we have for some $C>0$ \ce {\bf E}|X_t|^2\leq f(t)\leq
C\left[1+t^{\frac{1}{1-r/2}}\right], \de where the right side is
independent of $x_0$.

Since $p_t$ is irreducible and strong Feller, we have for any $b,
a>0$ and $t>0$ \ce \inf\limits_{x_0\in
B(0,b)}p_t\left(x_0,B(0,a)\right)>0. \de The second conclusion then
follows from \cite[Theorem 2.5 (b) and Theorem 2.7]{gm}.

\textbf{Acknowledgements:}

The author is very grateful to Professor Xicheng Zhang for his
valuable discussions and also wish to thank referees for suggestions and improvements.


\begin{thebibliography}{999}

\bibitem{c} S. Cerrai: Second Order PDE's in Finite and Infinite dimension. A Probabilistic
Approach, Lecture Notes in Math., vol. 1762, Springer-Verlag,
Berlin, 2001, x+330pp.

\bibitem{dz} G. Da Prato, J. Zabczyk: Ergodicity for Infinite Dimensional Systems, Cambridge
University Press, 1996.

\bibitem{gm} B. Goldys, B. Maslowski: Exponential ergodicity for stochastic reaction-
diffusion equations, Stochastic partial differential equations and
applications-VII, 115-131, Lect. Notes Pure Appl. Math., 245,
Chapman Hall/CRC, Boca Raton, FL, 2006.

\bibitem{hj} R. A. Horn, C. R. Johnson: Matrix Analysis, Cambridge University Press, 1985.

\bibitem{iw} N. Ikeda, S. Watanabe: Stochastic differential equations
and diffusion processes, 2nd ed., North-Holland/Kodanska,
Amsterdam/Tokyo, 1989.

\bibitem{Ku} A. M. Kulik: Exponential ergodicity of the solutions to SDE's with
a jump noise, {\it Stochastic Process. Appl}, 119(2009)602-632.

\bibitem{kw} Youngmff Kwon and Chaniio I. Ff: Strong feller property and irreducibility of
diffusions with jumps, {\it Stochastics An international Journal of
Probability and Stochastic Processes,} 67(1999)147-157.

\bibitem{ma1} H. Masuda: Ergodicity and exponential $\beta$-mixing bounds for multidimensional
diffusions with jumps, {\it Stochastic Processes and their
Applications,} 117(2007)35-56.

\bibitem{ma2} H. Masuda: Erratum to:``Ergodicity and exponential $\beta$-mixing bounds for multidimensional
diffusions with jumps"[Stochastic Process. Appl. 117(2007)35-56],
{\it Stochastic Processes and their Applications,} 119(2009)676-678.

\bibitem{pw} E. Priola and F. Wang: Gradient estimates for diffusion
semigroups with singular coefficients, {\it Journal of Functional
Analysis}, 236(2006)244-264.

\bibitem{rwz} J. Ren, J. Wu and X. Zhang: Exponential ergodicity of non-Lipschitz multivalued
stochastic differential equations, {\it Bull. Sci. Math.},
134(2010)391-404.

\bibitem{si} R. Situ: Theory of stochastic differential equations with jumps and applications,
Springer, 2005.

\bibitem{wa} J. Wang: Regularity of semigroups generated by L\'evy
type operators via coupling, {\it Stochastic Process. Appl.},
120(2010)1680-1700.

\bibitem{z1} X. Zhang: Exponential ergodicity of non-Lipschitz stochastic differential equations,
{\it Proceedings of the American Mathematical Society,}
137(2009)329-337.

\end{thebibliography}
\end{document}